\documentclass[12pt]{conm-p-l}
\usepackage{amscd}

\newcommand{\lab}[1]{\label{#1}}

\newtheorem{Thm}{Theorem}[section]
\newtheorem{Prop}[Thm]{Proposition}
\newtheorem{Lem}[Thm]{Lemma}
\newtheorem{Cor}[Thm]{Corollary}

\theoremstyle{remark}
\newtheorem{Rem}[Thm]{Remark}
\newtheorem*{Ack}{Acknowledgment}

\theoremstyle{definition}

\newtheorem{Def}[Thm]{Definition}

\newtheorem{Cond}[Thm]{Condition}

\newcommand{\bbR}{{\mathbb{R}}}

\newcommand{\cs}{{\rm CS}}

\newcommand{\de}{\partial}

\newcommand{\heta}{\hat\eta}

\newcommand{\Tr}{\operatorname{Tr}}

\newcommand{\wcC}{\widetilde{\mathcal{C}}}

\newcommand{\xinfty}{{x_\infty}}

\newcommand{\calS}{\mathcal{S}}
\newcommand{\calC}{\mathcal{C}}

\newcommand{\ord}{\operatorname{ord}}
\newcommand{\sln}{\operatorname{sln}}

\newcommand{\piHarm}{\operatorname{\pi_{\rm Harm}}}

\newcommand\qq{\rm}
\newcommand\cmp[1]{{\qq Commun.\ Math.\ Phys.\ \bf #1}}
\newcommand\jmp[1]{{\qq J.\ Math.\ Phys.\ \bf #1}}

\newcommand\anm[1]{{\qq Ann.\ Math.\ \bf #1}}

\newcommand\jdg[1]{{\qq J.\ Diff.\ Geom.\ \bf #1}}

\begin{document}
\title[Configuration Space Integrals and Invariants]
{Configuration Space Integrals and Invariants for 3-Manifolds
and Knots}
\author{Alberto S.~Cattaneo}
\address{Dipartimento di Matematica ``F.~Enriques" --- 
Universit\`a degli Studi di Milano ---
via Saldini, 50 --- I-20133 Milano --- Italy}
\email{cattaneo@elanor.mat.unimi.it}
\thanks{Supported by MURST and partially by INFN}



\maketitle

\section{Introduction}
In this paper we give a brief description of the way 
proposed in \cite{BC}
of associating invariants of both 3-dimensional
rational homology spheres (r.h.s.)\ and knots
in r.h.s.'s
to certain combinations of trivalent diagrams. In addition,
we discuss the relation between this construction and Kontsevich's
proposal \cite{K}.

The same diagrams appear in the LMO invariant \cite{LMO} for
3-manifolds, and it would be very interesting to know if there exists
any relationship between the two approaches in the case of r.h.s.'s.

The reason for restricting here to r.h.s.'s is
quite technical, as will be clear in Sect.~\ref{sec-irhs}.
Until that point,
without any loss of generality we can assume $M$ to be
any connected, compact, closed, oriented 3-manifold. 

Our construction
yields the invariants in terms of integrals over a suitable 
compactification of the configuration space of points on $M$.
More precisely, the
number of points corresponds to the number of vertices in the 
trivalent diagram, and the integrand is obtained by associating to
each edge in the diagram a certain 2-form that represents the
integral kernel of an ``inverse of the exterior derivative $d$."

One reason for constructing invariants in terms of ``$d^{-1}$" comes
from perturbative Witten--Chern--Simons theory \cite{W}. 
(More precisely,
one should invert the covariant derivative with respect to a flat
connection; so the present construction is related to the 
trivial-connection contribution.)

Another reason, which is perhaps more transparent to topologists, relies
on the definition of the linking number of two curves in $\bbR^3$
as the intersection number of one curve with any surface cobounding
the other. Thus, linking number may be defined in terms of an inverse
of the boundary operator $\de$. If one wants to represent the linking
number with an integral formula (viz., Gauss's formula), then 
one must consider the Poincar\'e duals of curves and apply ``$d^{-1}$."

The exterior derivative $d$ is neither injective nor surjective.
Thus, to invert it, one must restrict it to the complement of its kernel
and invert it on its image. Notice that one needs an explicit choice
of the complement of the kernel, and this introduces an element of
arbitrariness in the construction. Actually, our main task
will be to prove that the invariants we define are really independent 
of this arbitrary choice.

A general way of defining the inverse
of $d$ is by introducing a {\em parametrix}, i.e., a linear operator
on $\Omega^*(M)$ that decreases by one the form degree and satisfies
the following equation:
\begin{equation}
d\circ P + P\circ d = I - S,
\lab{dP}
\end{equation}
where $I$ is the identity operator and $S$ is a suitable smooth
operator such that \eqref{dP} has a solution. 

The definition
of the parametrix is far from unique. For the choice of $S$ is to a large
extent arbitrary and, moreover, given a solution $P$, 
we get another solution
of the form $P+d\circ Q - Q\circ d$ for any linear operator $Q$ that
decreases by two the form degree. These ambiguities reflect the ambiguities
in defining an inverse of $d$.

\begin{Rem}
One possible choice for
the parametrix is the Riemannian parametrix $P_g$, which is based on
the choice of a Riemannian metric $g$:
\begin{equation}
P_g = d^*\circ (\Delta+\piHarm)^{-1},
\qquad\alpha\in\Omega^*(M).
\lab{Pg}
\end{equation}
where $*$ is the Hodge-$*$ operator, $\Delta = d^*d+dd^*$ is the Laplace
operator, and $\piHarm$ is the projection
to harmonic forms. $P_g$ satisfies \eqref{dP} with $S=\piHarm$. 
\end{Rem}

The main point is now to define an integral kernel (of course not unique)
for $P$: i.e., we want to represent $P$ as a convolution. In the
language of differential topology, convolution can be written as
\[
P\alpha = -\pi_{2*}(\heta\ \pi_1^*\alpha),
\qquad\alpha\in\Omega^*(M),
\]
where $\pi_1$ and $\pi_2$ are the projections from $M\times M$
to either copy of $M$, and $\heta$ is the integral kernel for $P$
(wedges will be understood throughout).
By dimensional reasons it is clear that $\heta$ must be a 2-form on
$M\times M$. Because of the identity operator in \eqref{dP},
it is also clear that $\heta$ cannot be a smooth form. This is clearly
a problem since we want
to use $\heta$ to define smooth invariants for the manifold $M$. 
The solution consists in replacing $M\times M$ with a suitable
compactification $C_2(M)$ of the configuration space. 

The ambiguities in the definition of the parametrix imply that
the form $\heta$ is not unique. Thus, the main point will be
to prove that the invariants are independent of the choices involved
in the construction of $\heta$.
The plan of the construction is then the following:
\begin{enumerate}
\item One constructs a form $\heta$ to represent the integral kernel
of a parametrix with a prescription on its behavior on the boundary
of $C_2(M)$.
\item One introduces the unit interval $I$ as a space of parameters
to take care of the
arbitrary choices involved in this construction.
\item One associates to each trivalent diagram with $n$ vertices
a function on $I$ by integrating pullbacks of the form $\heta$
over a suitable compactification $C_n(M)$
(to be defined in Sect.~\ref{sec-cs}) of the configuration space
of $n$ points on $M$.
An invariant will then be a constant function on $I$.
\item
If the integrand form is closed, then the
differential of the corresponding function on $I$ is determined
only by contributions on the boundary of $C_n(M)$.
One shows then that if $M$ is a r.h.s.,
it is possible to associate a closed non-trivial integrand 
to each trivalent diagram.
\item Because of the prescribed behavior of $\heta$ on the boundary,
one can cancel the boundary contributions by summing up appropriate
combinations of diagrams.
\end{enumerate}
We discuss points 1.\ and 2.\ in Sect.~\ref{sec-constr} and
points 3., 4.\ and 5.\ in Sect.~\ref{sec-irhs}.

More concisely, the final statement is 
that to certain combinations (cocycles)
of trivalent diagrams it is possible to associate a
well-defined element of $H^{3n}(C_n(M),\de C_n(M))$, where
$M$ is a r.h.s.
The invariant for $M$ is then obtained by comparing this element
with the unit generator.

As we will see, the form $\heta$ is not closed, and this accounts for
the complications in point 4. This is due to the fact that the cohomology
of $M$ is not trivial. There are however cases when one can
obtain a closed form; viz.:
\begin{enumerate}
\item One can introduce a flat bundle $E$ over $M$
and consider the relevant covariant derivative instead of
the exterior derivative. If the bundle is non trivial, it may happen
that the complex $H^*(M;E)$ is acyclic. In this case, one can
construct a covariantly closed form to represent the inverse
of the covariant derivative.
\item If $M$ is a r.h.s., one can remove one point, thus obtaining a 
rational homology disc and, consequently, a closed form $\heta$.
\end{enumerate}

Case 1.\ was studied by Axelrod and Singer
\cite{AS} in the Riemannian framework of
eq.\ \eqref{Pg} (with $\Delta$ the covariant Laplace operator and,
by hypothesis, $\piHarm=0$). For a more general treatment of this case, 
s.\ \cite{BC2}. Notice, however, that this approach does not apply in general:
e.g., it does not even work for $S^3$ whose only flat connection is 
trivial.

Case 2.\ was proposed by Kontsevich \cite{K}. A realization of this
proposal for the simplest invariant---known as the 
$\Theta$-invariant---was then studied by Taubes \cite{T}. He also proved his 
version of the $\Theta$-invariant
to be trivial on integral homology spheres.
This rules out any relationship
with the LMO invariants which predict the $\Theta$-invariant 
on integral homology spheres
to be the Casson invariant.

In Sect.~\ref{sec-K} we will apply Kontsevich's proposal to the invariants
constructed in \cite{BC} (and Sect.~\ref{sec-irhs}), and will
compare our result with Taubes's.

\begin{Rem}
The above construction can be used to define invariants for
knots in a r.h.s.\ as well, s.\ \cite{BC} (and 
Sect.~\ref{sec-irhs}). 

Also in this case the construction is simplified if one gets
a closed form $\heta$. This happens, e.g., when $M=\bbR^3$ \cite{BT, AF},
or when $M=\Sigma\times [0,1]$ with $\Sigma$ 
a connected, compact,
closed, oriented 2-manifold \cite{AM}.

Another possibility is an approach \`a la Kontsevich when $M$ is
a r.h.s., s.\ Sect.~\ref{sec-K}.
\end{Rem}

\section{A compactification of configuration spaces}
\lab{sec-cs}

In this section we assume $M$ to be a $d$-dimensional
connected, compact, closed, oriented manifold.

The (open) configuration space $C_n^0(M)$ 
of $n$ points in $M$ is obtained by removing all diagonals
from the Cartesian product $M^n$. 

A $C^\infty$-compactification
for these spaces was proposed in \cite{AS} (generalizing the
algebraic compactification of \cite{FM}). This compactification is obtained
by taking the closure
\[
C_n(M) \doteq \overline{C_n^0(M)} \subset
M^n\times\prod_{\substack{S\subset\{1,2,\dots,n\}\\
|S|\ge2}} Bl(M^S,\Delta_S),
\]
where $Bl(M^S,\Delta_S)$ denotes the differential-geometric
blowup obtained by replacing the principal diagonal $\Delta_S$
in $M^S$ with its unit normal bundle $N(\Delta_S)/\bbR^+$.
Notice that $\Delta_S$ is diffeomorphic to $M$ and that 
\[N(\Delta_S)
\simeq TM^{\oplus S}/\text{global translations}.
\] 
Thus,
the boundary of $Bl(M^S,\Delta_S)$ is a bundle over $M$
associated to the tangent bundle.

Because of all the
blowups, the spaces $C_n(M)$ turn out to be manifolds with corners
($C_2(M)$ is simply a manifold with boundary). The codimension-one
components of the boundary of $C_n(M)$ are labeled by subsets
of $\{1,\dots,n\}$. By permuting the factors, we can always
put ourselves in the case when this subset is $\{1,\dots,k\}$,
$2\le k\le n$. We will denote by $\calS_{n,k}\subset \de C_n(M)$
a face of this kind. Then we have the following functorial description
of $\calS_{n,k}$:
\begin{equation}
\begin{CD}
\calS_{n,k} \simeq (\pi_1)^{-1}\widehat{C}_k(TM) @>>> \widehat{C}_k(TM)\\
 @VVV  @VVV\\
 C_{n-k+1}(M) @>>{\pi_1}> M
\end{CD}\lab{calS}
\end{equation}
Here $\pi_1$ is the projection onto the first copy of $M$
(i.e., where the first $k$ points have collapsed)
and $\widehat{C}_k(TM)$ is a bundle associated to the tangent
bundle of $M$ whose fiber $\widehat{C}_k(\bbR^d)$, $d=\dim M$,
is obtained from $(\bbR^d)^k/G$---$G$ being the group
of global translations and scalings---by 
blowing up all diagonals. More precisely, 
\[
\widehat{C}_k(\bbR^d) = \overline{{C_k^0(\bbR^d)}
/G }
\subset\left(
(\bbR^d)^k\times\prod_{\substack{S\subset\{1,2,\dots,k\}\\
|S|\ge2}} Bl((\bbR^d)^S,\Delta_S)\right)/G.
\]
Notice that each $\widehat{C}_k(\bbR^d)$ is a compact manifold
with corners. In the simplest case we have $\widehat{C}_2(\bbR^d)
=S^{d-1}$. 

The diagonal action of $SO(d)$
on $(\bbR^d)^k$ descends to $\widehat{C}_k(\bbR^d)$. If
we choose a Riemannian metric on $M$, we can write
\[
\widehat{C}_k(TM) =
OM \times_{SO(d)}\widehat{C}_k(\bbR^d),
\]
where $OM$ is the orthonormal frame bundle of $TM$.
In particular, when $n=2$, we have
\[
C_2(M) = Bl(M\times M,\Delta),
\]
and $\de C_2(M) \simeq S(TM) = OM\times_{SO(d)}S^{d-1}$.

\section{The construction of a parametrix}
\lab{sec-constr}

In this section we assume that $M$ is
any connected, compact, closed, oriented 3-manifold.

We start considering the
following commutative diagram:
\[
\begin{CD}
\de C_2(M) @>{\iota^\de}>> C_2(M) \\
@V{\pi^\de}VV @VV{\pi}V \\
\Delta @>>{\iota^\Delta}> M\times M
\end{CD}
\]
Then we define the involution $T$ that exchanges the
factors in $M\times M$. By abuse of notation, we will denote
by $T$ also the corresponding involution on $C_2(M)$ and on its
boundary $\de C_2(M)$. On the latter $T$ acts as the antipodal
map on the fiber crossed with the identity on the base. We will
denote by $H^*_\pm$ the $+$ and $-$ eigenspaces of $T$ in the 
cohomology of any of the above spaces.

We will denote by $\chi_\Delta\in\Omega^3(M\times M)$ 
a representative of the Poincar\'e dual of the diagonal $\Delta$.
Since $[\chi_\Delta]\in H^3_-(M\times M)$, there is really no
loss of generality in choosing an odd representative.

On the sphere bundle $\de C_2(M)\to\Delta$, one can introduce
a global angular form $\eta$; i.e., a form $\eta\in
\Omega^2(\de C_2(M))$ with the following properties:
\begin{enumerate}
\item the restriction of $\eta$ to each fiber is a generator
of the cohomology of the fiber;
\item $d\eta=-\pi^{\de*}e$, where $e$ is a representative of
the Euler class.
\end{enumerate}
Since $M$ is 3-dimensional, the Euler class is trivial.
Moreover, since $H^2_+(S^2)=0$, we may choose the global angular
form to be odd. 
Since $T$ acts as the identity on the base, 2.\ is then replaced by
$d\eta=0$. We have the following
\begin{Prop}
Given an odd global angular form $\eta$ and an odd representative
$\chi_\Delta$ of the Poincar\'e dual of the diagonal $\Delta$
in $M\times M$,
there exists a form $\heta\in\Omega^2(C_2(M))$ with the following
properties: \begin{subequations}\lab{propheta}
\begin{align}
d\heta &= \pi^*\chi_\Delta,\\
\iota_\de^*\heta &= -\eta,\\
T^*\heta &= -\heta.
\end{align}\end{subequations}
\lab{prop-heta}\end{Prop}
This is a simple generalization of the analogous proposition in \cite{BC}
for the case when $M$ is a r.h.s.
\begin{proof}
Let $U$ be a tubular neighborhood of $\Delta$ in $M\times M$,
and $\tilde U=\pi^{-1}U$ its preimage in $C_2(M)$. Then $\tilde U$
has the structure of $\de C_2(M)\times [0,1]$. Let us denote by
$\de_0\tilde U = \de C_2(M)$ and by $\de_1\tilde U = \pi^{-1}\de U$
the two boundary components of $\tilde U$.

Let $\rho$ be
a function on $\tilde U$ which is constant and equal to $-1$ in
a neighborhood of $\de_0\tilde U$ and is constant and equal
to $0$ in a neighborhood of $\de_1\tilde U$. Moreover, assume
that $\rho$ is even under the action of $T$.

Let $p:\tilde U\to \de C_2(M)$ be the natural projection. 
Then consider the form
$\tilde\eta = \rho\,p^*\eta$.
Since $\eta$ is a global angular form,
$d\tilde\eta = d\rho\,p^*\eta$ is a representative of the Thom
class of the normal bundle of $\Delta$. Therefore, if we extend
$\tilde\eta$ by zero on the whole of $C_2(M)$, we have that
$d\tilde\eta$ is the pullback of a representative of the Poincar\'e
dual of the diagonal.

This might not be our choice $\chi_\Delta$. 
Anyhow, $d\tilde\eta=\pi^*(\chi_\Delta + d\alpha)$,
and it is not difficult to check that one can choose 
$\alpha\in\Omega^2_-(M\times M)$. So we set
$\heta = \tilde\eta - \pi^*\alpha$,
and it is an immediate check that properties \eqref{propheta} hold.
\end{proof}

Notice that, as is clear from the proof, the definition of $\heta$
is not unique, even for fixed $\eta$ and $\chi_\Delta$.

With such a form $\heta$ we can finally define a parametrix. In fact,
let us denote by $\rho_1$ and $\rho_2$ the projections from $M\times M$
to each factor, and by $\pi_1$ and $\pi_2$ the 
corresponding projections from $C_2(M)$. By defining the
push-forward as acting from the left,
we have the following
\begin{Prop}
If $\heta\in\Omega^2(M)$ satisfies (\ref{propheta}a) and 
(\ref{propheta}b), then
\[
P\alpha = -\pi_{2*}(\heta\,\pi_1^*\alpha),
\qquad\alpha\in\Omega^*(M),
\]
is a parametrix with 
$S = -\rho_{2*}(\chi_\Delta\,\rho_1^*\alpha)$. 
\end{Prop}
The proof is a simple exercise
of fiber integration
in the case when the fiber has a boundary (s. \cite{BC}). 

\begin{Rem}
Property (\ref{propheta}c) is not needed in the definition of the
parametrix, but is natural and simplifies the writing of the 
invariants.
Notice incidentally that the 
{\em propagator}\/ in perturbative 
Witten--Chern--Simons theory is odd
under the action of $T$, being related to the {\em expectation
value}\/ of a connection 1-form placed at two different points.
\end{Rem}

\subsection{The global angular form}\lab{ssec-gaf}
To end with the construction, we have to specify a choice for the
global angular form.

Following \cite{AS} and \cite{BC}, we pick up
a Riemannian metric $g$. Then
$\de C_2(M) \simeq OM\times_{SO(3)} S^2$.
Let
\[
p:OM\times S^2 \to OM\times_{SO(3)} S^2
\]
be the natural projection,
and let $\theta$ be a connection form on $OM$ (i.e., a metric
connection).
By abuse of notation,
we will write $\theta$ also for its pullback to $OM\times S^2$.
\begin{Def}
We call a global angular form $\bar\eta\in\Omega^2(OM\times S^2)$
{\em equivariant}\/ if it is a polynomial in $\theta$ and 
$d\theta$---with coefficients in $\Omega^*(S^2)$---and
it is basic (i.e., $\bar\eta=p^*\eta$).
\end{Def}
For a given $\theta$, the equivariant global angular form is
unique and its $\theta$-independent part
is the $SO(3)$-invariant unit volume
form on $S^2$, which we will denote by $\omega$ throughout.
See \cite{BC} for an explicit construction.

\begin{Cond}
We assume that the restriction of $\heta$ to the boundary
is such that its pullback to $OM\times S^2$  
is the equivariant global angular form.\lab{cond-heta}
\end{Cond}

We denote by $\omega_{ij}$ the pullbacks of $\omega$ by
the projections $\widehat{C}_n(\bbR^3)\to\widehat{C}_2(\bbR^3)=S^2$.
Then a useful property proved in \cite{K} (s.\ also \cite{BC}) is
expressed by the following
\begin{Lem}[Kontsevich]
If $x_i$, $1\le i\le n$, is any coordinate in 
$\widehat{C}_n(\bbR^3)$, then, for any indices $j$ and $k$
($j\not=i$, $k\not=i$)
\[
\int_{x_i}\omega_{ij}\,\omega_{ik}=0.
\]
\lab{lem-K}
\end{Lem}

\subsection{Parametrizing the choices}
We have constructed an integral kernel $\heta$
for a parametrix depending on the following choices: 
a metric $g$, a compatible connection form $\theta$,
a representative $\chi_\Delta$
of the Poincar\'e dual of $\Delta$ (and a function $\rho$
and a 2-form $\alpha$ as in the proof of Prop.~\ref{prop-heta}).

To take care of these choices, we introduce the parameter space $I=[0,1]$.
Then, denoting by $\sigma$ the inclusion $M\times M\hookrightarrow
M\times M\times I$, we take $\chi_\Delta\in\Omega^3(M\times M\times I)$
such that $d\chi_\Delta=0$ and $\sigma^*\chi_\Delta$ is a representative
of the Poincar\'e dual of $\Delta$. (We treat similarly
$\rho$ and $\alpha$.)

As for $g$ and $\theta$, we operate as follows. We take a
block-diagonal
metric $g$ on $M\times I$
(i.e., a metric such that $g_{(m,t)}(v,w)=0$ for all
$v\in T_mM$ and $w\in T_tI$), and consider the orthonormal frame
bundle $\widetilde{OM}$ of $TM\to M\times I$. Then we choose
a connection form $\theta$ on $\widetilde{OM}$ and define
the equivariant global angular form on $\widetilde{OM}\times S^2$.

Using the projections $C_{n}(M)\times I\to C_2(M)\times I$,
we can pull back the form $\heta$ in $n(n-1)/2$ different ways
which we denote by $\heta_{ij}$.
\begin{Def}
We call a form on $\Omega^*(C_n(M))$ {\em special}\/ if it is a product
of pullbacks of $\heta$.
\end{Def}
Each special form can be graphically associated to a diagram,
each edge representing a pullback of $\heta$.

Let $\calS_{n,k}$ denote, as in the previous section, the face
in $\de C_n(M)$ corresponding to the collapse of the first $k$ 
points (we consider only this case since all other codimension-one faces
can be reduced to this one by simply applying a permutation of the
factors). Let $\pi^\calS$ denote the induced projection 
$\calS_{n,k}\times I\to C_{n-k+1}(M)\times I$, and let $\pi_1:
C_{n-k+1}(M)\times I\to M\times I$ denote the projection on the first
point (i.e., where the first $k$ points have collapsed).
Then we have the following
\begin{Lem}[Axelrod and Singer]
If $\alpha\in\Omega^*(\calS_{n,k})$ is the restriction of a special
form, then
\[
\pi^{\calS}_*\alpha = \beta\,\pi_1^*\gamma,
\]
where $\beta$ is special and $\gamma$ is either a constant or a
multiple of the first Pontrjagin form $p_1$ associated to $\theta$.
\lab{lem-AS}
\end{Lem}
For the proof, s.\ \cite{AS} or \cite{BC}. Using the same notations,
we also have the following
\begin{Cor}
$\gamma$ (and hence $\pi_1^*\gamma$)
is a constant
in the case when no parameter space is introduced.
\lab{cor-AS}
\end{Cor}

\section{An invariant for rational homology spheres}\lab{sec-irhs}
In this section we assume $M$ to be a 3-dimensional r.h.s. 
We then choose a representative $v$ of the unit generator of $H^3(M)$,
so we can take, as
the Poincar\'e dual of the diagonal
in $M\times M$,
\begin{equation}
\chi_\Delta = v_2 - v_1,
\lab{v21}
\end{equation}
where $v_i=\rho_i^*v$, and $\rho_i$, $i=1,2$, is the projection to
the $i$-th factor.

We now define our form $\heta$ as in the preceding
section---i.e., satisfying \eqref{propheta} and 
Condition~\ref{cond-heta}---with 
$\chi_\Delta$ as in \eqref{v21}. 

Next we consider the three projections $\pi_{ij}:C_3(M)\to C_2(M)$,
and write $\heta_{ij}=\pi_{ij}^*\heta$. Then \eqref{v21} implies
$d\heta_{ij} = v_j-v_i$.
Thus, we can define the following non-trivial closed form in 
$\Omega^2(C_3(M))$:
\begin{equation}
\heta_{123} \doteq \heta_{12}+\heta_{23}+\heta_{31}.
\end{equation}
Notice that on any configuration space $C_n(M)$, $n>2$, we can
analogously define closed forms $\heta_{ijk}$ for any triple
of distinct indices $i,j,k$.

Now consider graphs with numbered vertices, and 
set equivalent
to zero all graphs with an edge connecting a vertex to itself.
We have then an induced orientation of the edges (viz., each edge
is oriented from the lower to the higher end-point).

To each trivalent graph $\Gamma$ of the above type 
we can associate the following number:
\begin{equation}
A_\Gamma(M) \doteq \int_{C_{n+1}(M)} v_0\,
\prod_{(ij)\in E(\Gamma)}\heta_{ij0},
\lab{defAGamma}
\end{equation}
where $n$ is the number of vertices and $E(\Gamma)$ is the set of ordered
edges. The point labeled by $0$ is an extra point and not a vertex
of $\Gamma$.
We extend $A_\Gamma$  to combinations of graphs by linearity.

We are interested in the dependence of $A_\Gamma$ on
the choices in the construction of $\heta$. So we introduce a parameter
space $I$ as in the preceding section and consider $A_\Gamma$ as a
function on $I$. (As for $v$,
we take it in $\Omega^3(M\times I)$ and such that it
is closed and its restriction to $M$ is a representative of the
unit generator of $H^3(M)$.)
Then we consider the differential of this function. Since the integrand
form is closed, this differential is given only by boundary terms.
These are dealt with by using Lemmata \ref{lem-K} and \ref{lem-AS}.

\begin{Rem}
$A_\Gamma$ can be defined also if $M$ is not a r.h.s.
However, in this case the integrand form is not closed. So in 
differentiating $A_\Gamma$ we also have a bulk contribution which
we do not know how to deal with.
\end{Rem}

We now define a coboundary operator $\delta$ that acts on graphs
by contracting each edge one at a time, with a sign given by
the parity of the higher end-point. In \cite{BC}, it is shown that
$\delta$ is a coboundary operator. 

We call a {\em cocycle}\/ a $\delta$-closed combination of graphs.
We say that it is connected (trivalent) if all its terms are connected
(trivalent) graphs. 

Finally, we consider the Chern--Simons integral,
\[
\cs(M,f)=-\frac1{8\pi^2}\int_Mf^*\Tr\left({
\theta\,d\theta +\frac23\,\theta^3}\right),
\]
where $f$ is a section of $OM$ and $\theta$ is the same connection
form as in the construction of $\heta$. In \cite{BC}, the following
was proved:
\begin{Thm}
If\/ $\Gamma$ is a connected, trivalent cocycle, 
then there exists
a constant $\phi(\Gamma)$ such that
\[
I_\Gamma(M,f) = A_\Gamma(M) +\phi(\Gamma)\,\cs(M,f)
\]
is an invariant for the framed rational homology 3-sphere $M$.
\lab{thm-Gamma}
\end{Thm}

\begin{Rem}
Instead of defining
the equivariant global angular form,
one could repeat the previous construction by choosing a trivialization of 
$S(TM)$ and by defining the global angular form to be
the (pullback of the)
$SO(3)$-invariant unit volume form $\omega$ 
on $S^2$ (as suggested in \cite{K}). 

The equivariant treatment shows, 
cf.\ Thm.~\ref{thm-Gamma}, that under
a change of framing the invariants $A_\Gamma$ behave 
as multiples of the Chern--Simons integral.
In particular, they are invariant under a homotopic change of framing.
\lab{rem-triv}
\end{Rem}

All graphs in a cocycle have the same number $n$ of vertices.
If the cocycle is trivalent, this number is even. In this case,
one can define
\begin{equation}
\ord\Gamma = \frac n2.
\end{equation}
The constant $\phi(\Gamma)$ depends only on $\Gamma$ and not
on $M$. One can show \cite{AS,BC} that $\phi(\Gamma)=0$ if $\ord\Gamma$
is even. 
Moreover, in \cite{BC}
it is shown that $\phi(\Theta)=1/4$, with
\begin{equation}
A_\Theta = \int_{C_3(M)} v_0\,\heta_{012}^3.
\lab{defATheta}
\end{equation}
This allows for the definition
of the following unframed invariants (for $\Gamma\not=\Theta$):
\[
J_\Gamma(M) = A_\Gamma(M) - 4\,\phi(\Gamma)\, A_\Theta(M).
\]
A computation in \cite{BC} shows that $J_\Gamma(S^3)=J_\Gamma(SO(3))=0$
if $\ord\Gamma$ is odd.

\subsection{Knot invariants}
In \cite{BC}, invariants for knots in a r.h.s.\ are studied. 

If $K$ is an imbedding $S^1\hookrightarrow M$, one has induced imbeddings
$\widetilde C_n(S^1)\hookrightarrow C_n(M)$, where $\widetilde C_n(S^1)$
is the connected component of $C_n(S^1)$ defined by an
ordering of the points on $S^1$. The configuration space
$\widetilde\calC_{n,t}^K(M)$ of $n$ points on the knot and $t$ points
in $M$ is then defined by pulling back the bundle $C_{n+t}(M)\to
C_n(M)$.

All the forms introduced before can be pulled back to 
$\widetilde\calC_{n,t}^K$, and by abuse of notation
we will keep calling them with the same names. One should
keep in mind, however, that the pulled-back forms
depend on the imbedding $K$. This understood, one defines the
self-linking number
\[
\sln(K,M) \doteq \int_{\wcC_{2,0}^K(M)} \heta_{12}.
\]

Now consider graphs with a distinguished loop, which represents
the knot. We call {\em external}\/ the vertices and the edges
on this loop, and {\em internal}\/ all the others. To a trivalent
graph we can then associate the number
\begin{equation}
A_\Gamma(K,M) \doteq \int_{\wcC_{n,t+1}^K(M)} v_0\,
\prod_{(ij)\in I(\Gamma)} \heta_{ij0},
\lab{AKM}
\end{equation}
where $n$ and $t$ are the numbers of external and internal vertices
in $\Gamma$, and $I(\Gamma)$ is the set of internal
edges.
Again we extend \eqref{AKM} to combinations of graphs by linearity.

Next we define a coboundary operator $\delta$ as before with
the only difference that now $\delta$ does not contract internal
edges connecting two external vertices. Graph combinations killed
by $\delta$ are called cocycles. (An explicit computation of these
cocycles is presented in \cite{AF}.)

We name {\em prime}\/ a graph which is connected after removing
any pair of external edges (in \cite{BC} a graph of this kind
was called connected). 
A cocycle will be called prime (trivalent)
if all its terms are prime (trivalent).
In \cite{BC}, the following was proved:
\begin{Thm}
If $K$ is a knot in the rational homology 3-sphere $M$, and\/
$\Gamma$ is a prime, trivalent cocycle, 
then there exists a constant $\mu(\Gamma)$ such that
\[
I_\Gamma(K,M) = A_\Gamma(K,M) + \mu(\Gamma)\, \sln(K,M)
\]
is a knot invariant. Moreover, $\mu(\Gamma)=0$ if\/ $\ord\Gamma$ is
even.
\lab{thm-GammaK}
\end{Thm}

\section{Relationship with Kontsevich's proposal} \lab{sec-K}
As we recalled in the Introduction,
in \cite{K} Kontsevich proposed a different way of constructing
invariants for r.h.s.'s. His proposal differs from
our construction since $i)$ one point is removed from $M$ in order to 
make its rational homology trivial, 
and $ii)$ the global angular form on the boundary
of the configuration space is defined via a trivialization of $TM$.
Let us consider part $i)$ of the proposal first.

We introduce the compactified configuration space $C_n(M,\xinfty)$
of $n$ points on $M\backslash\xinfty$
(where $\xinfty$ is an arbitrary base point)
as the fiber of $C_{n+1}(M)\to M$ over
the point $\xinfty$ in the last copy of $M$; viz.:
\[
\begin{CD}
C_n(M,\xinfty) @>{p}>> C_{n+1}(M) \\
@VVV @VV{\pi_{n+1}}V \\
\xinfty @>>> M
\end{CD}
\]
Consider now the projections $\pi_{ij}:C_{n+1}(M)\to C_2(M)$,
$i<j\le n+1$, and 
set $p_{ij}=\pi_{ij}\circ p$ for $i<j\le n$ and
$p_{i\infty}=\pi_{i,n+1}\circ p$ for $i\le n$. Then we will
denote by $\heta_{ij}$ and $\heta_{i\infty}$ the pullbacks
of $\heta$ to $C_n(M,\xinfty)$.  Accordingly, we will define
$\heta_{ijk}$ and $\heta_{ij\infty}$. On the boundary faces
we will have the pullbacks of global angular forms 
$\eta_{ij}$ and $\eta_{i\infty}$. Observe that $\eta_{i\infty}=
\omega_{i\infty}$, with the notations of Lemma \ref{lem-K}.

We also have projections $C_n(M,\xinfty)\to C_n(M)$. The forms
$\heta_{ij}$ we have written above can also be seen as the pullbacks of the
forms with the same name on $C_n(M)$.

In particular, $C_1(M,\xinfty)$ is just $M$ blown up at $\xinfty$,
and $\de C_1(M,\xinfty)=S^2$. If we pull $v$ back by the projection
$\tau:C_1(M,\xinfty)\to M$, we get an exact form $\tau^*v=d w$, where
the two-form $w$ restricted to the boundary must be a representative
of the unit generator of $H^2(S^2)$. 
We may choose this representative to be $\omega$ 
since, by Thm.~\ref{thm-Gamma}, 
the explicit choice of $v$ does not affect the invariant.
We have then the following
\begin{Thm}
If\/ $\Gamma$ is a connected, trivalent cocycle 
and $M$ is a rational homology 3-sphere, then
\[
A_\Gamma(M) = A_\Gamma'(M) + B_\Gamma,
\]
with
\begin{align*}
A_\Gamma'(M) &= \int_{C_{n}(M,\xinfty)}
\prod_{(ij)\in E(\Gamma)}\heta_{ij\infty},\\
B_\Gamma &= \int_{\widehat C_{n+2}(\bbR^3)}
\omega_{0\infty}\, \prod_{(ij)\in E(\Gamma)}\omega_{ij0}.
\end{align*}
Moreover, if\/ $\ord\Gamma$ is odd, $B_\Gamma$ vanishes.
\lab{thm-AGamma'}
\end{Thm}

Notice that $B_\Gamma$ does not depend on $M$ or on any arbitrary
choice. Thus, 
even if it should not vanish,
it would just represent a constant shift in the invariant. 
As a consequence,
Thm.~\ref{thm-Gamma} holds with $A_\Gamma(M)$ replaced by
$A_\Gamma'(M)$.

\begin{proof}
We can pull back the integrand form in \eqref{defAGamma}
to $C_{n+1}(M,\xinfty)$ and integrate it over there.
Using the fact that the pullback of $v_0$ is exact, 
by Stokes's theorem we can rewrite \eqref{defAGamma}
as
\[
A_\Gamma(M)=\int_{\de C_{n+1}(M,\xinfty)} w_0\,
\prod_{(ij)\in E(\Gamma)}\heta_{ij0}.
\]
The codimension-one faces in $\de C_{n+1}(M,\xinfty)$ are
labeled by subsets of $\{0,1,\dots,n,\infty\}$. Denote by $\calS$
any of these subsets. 

Assume now that the cardinality of $\calS'=\calS\cap\{1,\dots,n\}$
is $k$.
Since points in $\calS'$ label vertices in the graph and the graph is 
trivalent, we have the relation
\[
3k=2e+e_0,
\]
where $e$ denotes the number of edges connecting points in $\calS'$,
and $e_0$ denotes the number of edges with exactly one
end-point in $\calS'$. Now we have four cases, 
according to 
\[
\calS\backslash\calS'= 
\begin{cases}
\{0\} &\text{(a)}\\
\emptyset &\text{(b)}\\
\{0,\infty\} &\text{(c)}\\
\{\infty\} &\text{(d)}\\
\end{cases}
\]
The cardinality $r$ of $\calS$ is then $k+1$, $k$, $k+2$ and $k+1$
respectively. The boundary face labeled by $\calS$ is
a bundle over $C_{n+2-r}(M,\xinfty)$ with projection $\pi^{\calS}$ and
fiber $\widehat{C}_r(\bbR^3)$.
So the fiber dimension is $3k-1$, $3k-4$, $3k+2$
and $3k-1$ respectively.

We now
write the integrand form $\alpha$ restricted to this boundary 
as $(\pi^{\calS*}\beta)\,\alpha'$, where $\beta$ is the product of
the pullbacks of $\heta$ corresponding to edges with at least
one end-point not in $\calS$, times $w_0$ in cases (a), (b) and (d).

In cases (b) and (d),
the term $\heta_{ij0}$ contributes
to $\alpha'$ only if both $i$ and $j$ are in $\calS'$. In case (a) and
(c), also terms with either $i$ or $j$ in $\calS'$ contribute.
Moreover, $w_0$  contributes to $\alpha'$ only in case (c). 
As a consequence,
the degree of $\alpha'$ will be: (a) $2e+2e_0$, (b) $2e$, 
(c) $2e+2e_0+2$, (d) $2e$.

By using all the above results, we see that the degree of $\gamma$
in Corollary \ref{cor-AS} is $e_0+1$, $4-e_0$, $e_0$ and $1-e_0$ 
respectively.
Since $\gamma$ must be a constant zero-form, we see that the contribution
of the face $\calS$ vanishes unless we are in case (b) with $e_0=4$,
in case (c) with $e_0=0$, or in case (d) with $e_0=1$. Notice, moreover,
that we can replace $\heta$ by $\omega$ in $\alpha'$. Thus, in the last
case above we conclude that the contribution vanishes 
by Lemma \ref{lem-K}. The first case is
taken care of by the fact that $\Gamma$ is a cocycle.

We are then left with case (c) and $e_0=0$. Since $\Gamma$ is connected,
there are only two possibilities: 1) only point $0$ has collapsed at
$\xinfty$, 2) all points have collapsed at $\xinfty$.
In case 1), $\alpha'=\omega_{0\infty}$ and the fiber is $S^2$. After this
trivial integration we get $A_\Gamma'(M)$.
Case 2) yields $B_\Gamma$. 

To prove that $B_\Gamma$ vanishes
if $\ord\Gamma = n/2$ is odd, consider the involution $x_i\to-x_i$,
$i=0,1,\dots,n,\infty$. All the pullbacks of $\omega$ change signs.
Since the number of edges is $3n/2$, the integrand form gets
the sign $(-1)^{3n/2+1}$. On the other hand, since
$\widehat C_{n+2}(\bbR^3)$ is $S^{3n+2}$ with some submanifolds
blown up, under the involution the orientation gets the sign
$(-1)^{n+1}$.
\end{proof}

In the particular case of the $\Theta$-invariant \eqref{defATheta},
we have
\[
A_\Theta'(M) = \int_{C_2(M,\xinfty)} \heta_{12\infty}^3.
\]
By our construction, $\heta_{12\infty}$ is a closed form on
$C_2(M,\xinfty)$ which reduces to the global angular form
when restricted to the faces $(1\infty)$, $(2\infty)$ and
$(12)$. 

As observed in remark \ref{rem-triv}, one can also
modify the construction by choosing a trivialization of
$S(TM)$ at the very beginning, and this corresponds to part $ii)$
of Kontsevich's proposal.
Invariance under homotopic
changes of framing is then guaranteed (while under non-homotopic
changes, the invariant behaves as $-1/4\,\cs$). 
In this case, we have the additional property that
$\heta_{12\infty}^2$ vanishes close to 
the faces $(1\infty)$, $(2\infty)$ and
$(12)$. However, close to $(12\infty)$ neither
$\heta_{12\infty}^2$ nor $\heta_{12\infty}^3$ vanish.

This is to be compared with Taubes's invariant
\[
\widetilde A_\Theta(M) \doteq \int_{C_2(M,\xinfty)} \omega^3,
\]
where $\omega$ is a 2-form on $C_2(M,\xinfty)$ with the following
properties:
\begin{enumerate}
\item $\omega$ restricted to the faces $(1\infty)$, $(2\infty)$ and
$(12)$ is a global angular form;
\item $\omega^2$ vanishes not only close to $(1\infty)$, $(2\infty)$ and
$(12)$ but also close to $(12\infty)$.
\end{enumerate}
The latter property is achieved only by choosing what Taubes names a
{\em singular framing}\/ for $T(M\backslash\xinfty)$.
As a consequence, $\omega^2$ (and hence $\omega^3$) is a form with
compact support, and Taubes's $\Theta$-invariant can actually
be defined as an integral over the uncompactified configuration
space $C_2^0(M,\xinfty)$.
Moreover, property 2.\ is crucial in Taubes's proof that his invariant
is trivial on integral homology spheres.

Now the main question is if there is any relationship between the
two different ways, $A'_\Theta$ and $\widetilde A_\Theta$,
of realizing Kontsevich's proposal for the $\Theta$-invariant.

\subsection{The case of knots}
Let us consider an imbedding $K$ of $S^1$
in the interior of $M\backslash\xinfty$.
This induces imbeddings $\widetilde C_n(S^1)\hookrightarrow C_n(M,\xinfty)$.
By pulling back the bundles $C_{n+t}(M,\xinfty)\to C_n(M,\xinfty)$,
we then obtain the configuration spaces $\calC_{n,t}^{K}(M,\xinfty)$.
We have the following
\begin{Thm}
If\/ $\Gamma$ is a prime, trivalent cocycle 
and $K$ is a knot in
the rational homology 3-sphere $M$, then
\begin{align*}
A_\Gamma(K,M) &= \int_{\wcC_{n,t}^{K}(M,\xinfty)} 
\prod_{(ij)\in I(\Gamma)} \heta_{ij\infty},\\
\sln(K,M) &= \int_{\wcC_{2,0}^{K}(M,\xinfty)} \heta_{12}.
\end{align*}
\end{Thm}

In particular, if $M=S^3$ and we choose the Euclidean metric
on $\bbR^3=S^3\backslash x_\infty$,
we recover Bott and Taubes's result \cite{BT}.
As a consequence, the anomaly coefficients $\mu(\Gamma)$ 
are the same in the two cases.

\begin{proof}
We work as in the proof of Thm.~\ref{thm-AGamma'}. The only difference
is that we must distinguish between the cases when the collapse is at a
point on $K$ or otherwise. 

Notice that, since $\xinfty$ does not belong to the image
of $K$, there is no such term as $B_\Gamma$. By the same reason,
when we consider a collapse at a point on $K$, we only have points
in $\{0,1,\dots,n+t\}$. If $0$ is involved, the term vanishes since
$w_0$ is basic and is a 2-form. If $0$ is not involved, reasoning as
in the proof of Thm.~\ref{thm-AGamma'} and applying Corollary
\ref{cor-AS}
shows that the term vanishes
unless $e_0=2$. But this is taken care of by the fact that $\Gamma$ is
a cocycle.
\end{proof}

\begin{Ack}
I thank J.~E.~Andersen,
R.~Bott, P.~Cotta-Ramusino, N.~Habegger, R.~Longoni, 
G.~Masbaum and P.~Vogel for very useful discussions.

For partial support at Madeira conference on ``Low Dimensional
Topology," January 1998,
I acknowledge the Center of Mathematical Sciences (CCM) PRAXIS XXI 
project. I especially thank H.~Nencka for her efforts in organizing
the conference.

I am finally thankful to the Institut de Mathematiques de Jussieu and to the
University of Nantes for warm hospitality and financial support.
\end{Ack}

\end{document}